\documentclass[a4paper,10pt]{amsart}
\usepackage[english]{babel}
\usepackage[applemac]{inputenc}
\usepackage{amsmath}
\usepackage{eucal}
\usepackage{amsthm}
\usepackage{amssymb}

\newtheorem{thm}{Theorem}
\newtheorem{prop}{Proposition}

\newtheorem{lem}{Lemma}
\theoremstyle{definition}
\newtheorem{dfn}{Definition}
\theoremstyle{remark}
\newtheorem{rem}{Remark}

\DeclareMathOperator{\rank}{\text{\rm rank}}
\DeclareMathOperator{\codim}{\text{\rm codim}}
\DeclareMathOperator{\length}{\text{\rm length}}

\author{Simone Diverio}
\email{sdiverio@fourier.ujf-grenoble.fr}
\address{Institut Fourier\\
Université de Grenoble I\\
BP 74, F-38402, Saint Martin d'Hères\\
France.}
\curraddr{Istituto \lq\lq Guido Castelnuovo\rq\rq{}\\
Università di Roma \lq\lq La Sapienza\rq\rq{}\\
P.le Aldo Moro, 2,  00185 Roma\\
Italia.}

\title[DIFFERENTIAL EQUATIONS ON COMPLEX PROJECTIVE HYPERSURFACES]{DIFFERENTIAL EQUATIONS ON COMPLEX PROJECTIVE HYPERSURFACES OF LOW DIMENSION}

\begin{document}

\begin{abstract}
Let $n=2,3,4,5$ and let $X$ be a smooth complex projective hypersurface of $\mathbb P^{n+1}$. In this paper we find an effective lower bound for the degree of $X$, such that every holomorphic entire curve in $X$ must satisfy an algebraic differential equation of order $k=n=\dim X$, and also similar bounds for order $k>n$. Moreover, for every integer $n\ge 2$, we show that there are no such algebraic differential equations of order $k<n$ for a smooth hypersurface in $\mathbb P^{n+1}$.
\end{abstract}

\maketitle

\section{Introduction}
Let $X\subset\mathbb P^{n+1}$ be a complex projective hypersurface, with $\deg X=d$. In 1970 S. Kobayashi conjectured \cite{Kobayashi70} that $X$ is hyperbolic provided $X$ is generic and $d$ is large enough. During recent years, several efforts have been made to treat both the low-dimension cases (with a special attention to the lower bound for the degree) and the general one.

For instance, \cite{D-EG00} prove the conjecture for very generic surfaces in $\mathbb P^3$ of degree greater or equal to $21$, \cite{Rousseau05} proves a weaker form (namely weak analytic hyperbolicity) for generic hypersurfaces in $\mathbb P^4$ of degree greater than or equal to $593$ and \cite{Siu04} announces a complete proof of the Kobayashi conjecture in any dimension, for $d\gg n$.

The major technique in the study of the Kobayashi hyperbolicity question is probably to analyse the existence of sections of the bundle $E_{k,m}^{\rm GG}T^*_X$ of differentials of order $k$ and weighted degree $m$, as defined in \cite{G-G79}. This technique was later refined by J.-P. Demailly in \cite{Demailly95}, with the introduction of the bundle $E_{k,m}T^*_X$ of invariant differential operators.

Let $X$ be a compact complex manifold. Then $X$ is Kobayashi hyperbolic if and only if there are no non-constant entire holomorphic curves in $X$ (Brody's criterion). The general philosophy is that global holomorphic sections of $E_{k,m}T^*_X$ vanishing on a fixed ample divisor give rise to global algebraic differential equations that every entire holomorphic curve must satisfy.

It is known by \cite{D-EG00} that every smooth surface in $\mathbb P^3$ of degree greater or equal to $15$ has such differential equations of order two. For the dimension three case \cite{Rousseau06b} observed that one needs to look for order three equations since one has in general the vanishing of symmetric differentials and invariant $2$-jet differentials for smooth hypersurfaces in projective $4$-space. On the other hand \cite{Rousseau06b} shows the existence of global invariant $3$-jet differentials vanishing on an ample divisor on every smooth hypersurface $X$ in $\mathbb P^4$, provided that $\deg X\ge 97$. 

The existence of these global sections is shown by means of a delicate algebraic study of the bundle $E_{3,m}T^*_X$ conducted in \cite{Rousseau06a}, which permits to compute the Euler characteristic. Then by a laborious estimate of the dimension of the higher cohomology groups, one can find a positive lower bound for the dimension of the space of global sections (at least for $m$ large).

In the present paper, we first of all generalise to all dimension the non-existence of global sections of invariant jet differentials of order less than the dimension of the ambient variety. In fact, we prove the following

\begin{thm} 
Let $X\subset\mathbb P^{n+1}$ be a smooth hypersurface. Then 
$$
H^0(X,E_{k,m}T^*_X)=0
$$
for all $m\ge 1$ and $1\le k\le n-1$. In other words, on a smooth projective hypersurface there are no global invariant jet differentials of order less than its dimension.
\end{thm} 

The idea of the proof is to exclude, in the direct sum decomposition into irreducible Gl$(T^*_X)$-representation of the graded bundle $\operatorname{Gr}^\bullet E_{k,m}T^*_X$, the existence of Schur powers $\Gamma^{(\lambda_1,\dots,\lambda_n)}T^*_X$, with $\lambda_n>0$, and then to use a vanishing theorem due to P. Br\"uckmann and H.-G. Rackwitz for Schur powers of the cotangent bundle of smooth projective complete intersections.

On the other hand, in the direction of existence of global invariant jet differentials, we get a slightly better bound for smooth hypersurfaces in $\mathbb P^4$ and a new result for smooth hypersurfaces in $\mathbb P^5$ and $\mathbb P^6$.

\begin{thm}
Let $X\subset\mathbb P^{n+1}$, $n=2,3,4,5$, be a smooth hypersurface of degree $d$ and $A\to X$ an ample line bundle. Then
$$
H^0(X,E_{k,m}T^*_X\otimes A^{-1})\ne 0
$$
for $m$ large enough and $k\ge n$ if $d$ is sufficiently big, and therefore every holomorphic entire curve $f\colon\mathbb C\to X$ must satisfy the corresponding algebraic differential equation. 

\noindent
Moreover, we have the following effective lower bounds for the degree $d$ (depending on the values of $n$ and~$k$ which are the entries of the table):
$$ 
\begin{array}{|c||c|c|c|c|c|}
\hline
\bold{n} \backslash \bold{k} & 1 & 2 & 3 & 4 & 5  \\ 
\hline 
\hline
2  & - & 18 & 16 & 16 & 16 \\
\hline
3 & - & - & 82 & 74 & 74 \\
\hline
4 & - & - & - & 329 & 298\\
\hline
5 & - & - & - & - & 1222\\
\hline
\end{array} 
$$
\end{thm}

Here the proof is achieved without any Euler characteristic computation, thanks to the algebraic version \cite{Trapani95} of Demailly's holomorphic Morse inequalities. In our case, the Morse inequalities are applied to a particular subalgebra of $E_{k,m}T^*_X$ which is easier to compute than the full algebra $\bigoplus_m E_{k,m}T^*_X$ itself~-- although its positivity properties are probably not as good.

We would like to point out that in this way no higher cohomology computations nor any algebraic study of the jet bundle are needed. However, even if {\sl a priori} these techniques should work in higher dimension, the amount of computations needed to get the result blows-up rapidly when the dimension increases.

\subsection{Acknowledgements}
I would like to thank my two thesis directors, Prof. Jean-Pierre Demailly and Prof. Stefano Trapani for their useful help, encouragement and their extreme patience and very nice attitude. And of course for all the formal talk and informal chats we had which introduced me in this subject.

Thanks also to Andrea Maffei for his \lq\lq algebraic support\rq\rq{} and to Erwan Rousseau for having generously shared with me some of his ideas about the vanishing of jet differentials.

\section{Background Material and Preliminaries}

In this section, we follow very closely \cite{Demailly95}.

\subsection{Jet Differentials}

Let $(X,V)$ be a directed manifold, {\sl i.e.} a pair where $X$ is a complex manifold and $V\subset T_X$ a holomorphic subbundle (non necessarily integrable) of the tangent bundle. The bundle $J_kV$ is the bundle of $k$-jets of holomorphic curves $f\colon(\mathbb C,0)\to X$ which are tangent to $V$, i.e., such that $f'(t)\in V_{f(t)}$ for all $t$ in a neighbourhood of $0$, together with the projection map $f\mapsto f(0)$ onto $X$.

Let $\mathbb G_k$ be the group of germs of $k$-jets of biholomorphisms of $(\mathbb C,0)$, that is, the group of germs of biholomorphic maps
$$
t\mapsto\varphi(t)=a_1\, t+ a_2\, t^2+\cdots+a_k\,t^k,\quad a_1\in\mathbb C^*,a_j\in\mathbb C,j\ge 2,
$$
in which the composition law is taken modulo terms $t^j$ of degree $j>k$. Then $\mathbb G_k$ admits a natural fiberwise right action on $J_kV$ consisting of reparametrizing $k$-jets of curves by a biholomorphic change of parameter. Moreover the subgroup $\mathbb H\simeq\mathbb C^*$ of homotheties $\varphi(t)=\lambda\,t$ is a (non normal) subgroup of $\mathbb G_k$ and we have a semidirect decomposition $\mathbb G_k=\mathbb G'_k\ltimes\mathbb H$, where $\mathbb G'_k$ is the group of $k$-jets of biholomorphisms tangent to the identity. The corresponding action on $k$-jets is described in coordinates by
$$
\lambda\cdot(f',f'',\dots,f^{(k)})=(\lambda f',\lambda^2 f'',\dots,\lambda^k f^{(k)}).
$$
As in \cite{G-G79}, we introduce the vector bundle $E_{k,m}^{GG}V^*\to X$ whose fibres are complex valued polynomials $Q(f',f'',\dots,f^{(k)})$ on the fibres of $J_kV$, of weighted degree $m$ with respect to the $\mathbb C^*$ action defined by $\mathbb H$, that is, such that
$$
Q(\lambda f',\lambda^2 f'',\dots,\lambda^k f^{(k)})=\lambda^mQ(f',f'',\dots,f^{(k)}),
$$
for all $\lambda\in\mathbb C^*$ and $(f',f'',\dots,f^{(k)})\in J_kV$.

Next, we define the bundle of Demailly-Semple jet differentials (or invariant jet differentials) as a subbundle of the Green-Griffiths one.

\begin{dfn}[\cite{Demailly95}]
The \emph{bundle of invariant jet differentials of order $k$ and degree $m$} is the subbundle $E_{k,m}V^*\subset E_{k,m}^{GG}V^*$ of polynomial differential operators $Q(f',f'',\dots,f^{(k)})$ which are invariant under arbitrary changes of parametrization, i.e., for every $\varphi\in\mathbb G_k$
$$
Q((f\circ\varphi)',(f\circ\varphi)'',\dots,(f\circ\varphi)^{(k)})=\varphi'(0)^m\,Q(f',f'',\dots,f^{(k)}).
$$
Alternatively, $E_{k,m}V^*=(E_{k,m}^{GG}V^*)^{\mathbb G'_k}$ is the set of invariants of $E_{k,m}^{GG}V^*$ under the action of $\mathbb G'_k$.
\end{dfn}

We now define a filtration on $E_{k,m}^{GG}V^*$. A coordinate change $f\mapsto\Psi\circ f$ transforms every monomial $(f^{(\bullet)})^\ell=(f')^{\ell_1}(f'')^{\ell_2}\cdots(f^{(k)})^{\ell_k}$ of partial weighted degree $|\ell|_s:=\ell_1+2\ell_2+\cdots+s\ell_s$, $1\le s\le k$, into a polynomial $((\Psi\circ f)^{(\bullet)})^\ell$ in $(f',f'',\dots,f^{(k)})$, which has the same partial weighted degree of order $s$ if $\ell_{s+1}=\cdots=\ell_k=0$ and a larger or equal partial degree of order $s$ otherwise. Hence, for each $s=1,\dots,k$, we get a well defined decreasing filtration $F_s^\bullet$ on $E_{k,m}^{GG}V^*$ as follows:
$$
F_s^p(E_{k,m}^{GG}V^*)=
\left\{
\begin{matrix}
\text{$Q(f',f'',\dots,f^{(k)})\in E_{k,m}^{GG}V^*$ involving} \\
\text{only monomials $(f^{(\bullet)})^\ell$ with $|\ell|_s\ge p$}
\end{matrix}
\right\},
\quad\forall p\in\mathbb N.
$$
The graded terms $\operatorname{Gr}^p_{k-1}(E_{k,m}^{GG}V^*)$, associated with the filtration $F_{k-1}^p(E_{k,m}^{GG}V^*)$, are precisely the homogeneous polynomials $Q(f',f'',\dots,f^{(k)})$ whose all monomials $(f^{(\bullet)})^\ell$ have partial weighted degree $|\ell|_{k-1}=p$; hence, their degree $\ell_k$ in $f^{(k)}$ is such that $m-p=k\ell_k$ and $\operatorname{Gr}^p_{k-1}(E_{k,m}^{GG}V^*)=0$ unless $k|m-p$. Looking at the transition automorphisms of the graded bundle induced by the coordinate change $f\mapsto\Psi\circ f$, it turns out that $f^{(k)}$ behaves as an element of $V\subset T_X$ and, as a simple computation shows, we find
$$
\operatorname{Gr}^{m-k\ell_k}_{k-1}(E_{k,m}^{GG}V^*)=E_{k-1,m-k\ell_k}^{GG}V^*\otimes S^{\ell_k}V^*.
$$
Combining all filtrations $F^\bullet_s$ together, we find inductively a filtration $F^\bullet$ on $E_{k,m}^{GG}V^*$ such that the graded terms are
$$
\operatorname{Gr}^{\ell}(E_{k,m}^{GG}V^*)=S^{\ell_1}V^*\otimes S^{\ell_2}V^*\otimes\cdots\otimes S^{\ell_k}V^*,\quad\ell\in\mathbb N^k,|\ell|_k=m.
$$
Moreover there are natural induced filtrations $F^p_s(E_{k,m}V^*)=E_{k,m}V^*\cap F_{s}^p(E_{k,m}^{GG}V^*)$ in such a way that
$$
\operatorname{Gr}^\bullet(E_{k,m}V^*)=\left(\bigoplus_{|\ell|_k=m}S^{\ell_1}V^*\otimes S^{\ell_2}V^*\otimes\cdots\otimes S^{\ell_k}V^*\right)^{\mathbb G'_k}.
$$
Here, we remark that, in general, it is a major unsolved problem to find the decomposition of $\operatorname{Gr}^\bullet(E_{k,m}V^*)$ into irreducible Gl$(V^*)$-representations. This is easy for $k\le 2$ (since the summands in the graded bundle do not mix up under the action of $\mathbb G'_k$); \cite{Rousseau06a} also found the formula for $k=\dim X=3$ and $V=T_X$, thanks to a deeper study based on a theorem of V.~Popov on invariant theory.

\subsection{Projectivized $\bold k$-Jet Bundles}

Here we explain the construction of the tower of projectivized bundles which provides a (relative) smooth compactification of $J^{\text{reg}}_kV/\mathbb G_k$, where $J^{\text{reg}}_kV$ is the bundle of regular $k$-jets tangent to $V$, that is, $k$-jets such that $f'(0)\ne 0$.

Let $(X,V)$ be a directed manifold, with $\dim X=n$ and $\rank V=r$. With $(X,V)$ we associate another directed manifold $(\widetilde X,\widetilde V)$ where $\widetilde X=P(V)$ is the projectivized bundle of lines of $V$, $\pi\colon\widetilde X\to X$ is the natural projection and $\widetilde V$ is the subbundle of $T_{\widetilde X}$ defined fiberwise as
$$
\widetilde V_{(x_0,[v_0])}\overset{\text{def}}=\{\xi\in T_{\widetilde X,(x_0,[v_0])}\mid\pi_*\xi\in\mathbb C.v_0\},
$$
$x_0\in X$ and $v_0\in T_{X,x_0}\setminus\{0\}$. 
We also have a \lq\lq lifting\rq\rq{} operator which assigns to a germ of holomorphic curve $f\colon(\mathbb C,0)\to X$ tangent to $V$ a germ of holomorphic curve $\widetilde f\colon(\mathbb C,0)\to\widetilde X$ tangent to $\widetilde V$ in such a way that $\widetilde f(t)=(f(t),[f'(t)])$.

To construct the projectivized $k$-jet bundle we simply set inductively $(X_0,V_0)=(X,V)$ and $(X_k,V_k)=(\widetilde X_{k-1},\widetilde V_{k-1})$. Of course, we have for each $k>0$ a tautological line bundle $\mathcal O_{X_k}(-1)\to X_k$ and a natural projection $\pi_k\colon X_k\to X_{k-1}$. We call $\pi_{j,k}$ the composition of the projections $\pi_{j+1}\circ\cdots\circ\pi_{k}$, so that the total projection is given by $\pi_{0,k}\colon X_k\to X$.
We have again for each $k>0$ short exact sequences
\begin{equation}\label{ses1}
0\to T_{X_k/X_{k-1}}\to V_k\to\mathcal O_{X_k}(-1)\to 0,
\end{equation}
\begin{equation}\label{ses2}
0\to\mathcal O_{X_k}\to\pi_k^*V_{k-1}\otimes\mathcal O_{X_k}(1)\to T_{X_k/X_{k-1}}\to 0
\end{equation}
and $\rank V_k=r$, $\dim X_k=n+k(r-1)$. Here, we also have an inductively defined $k$-lifting for germs of holomorphic curves such that $f_{[k]}\colon(\mathbb C,0)\to X_k$ is obtained as $f_{[k]}=\widetilde f_{[k-1]}$.

The following theorem justifies in some sense the construction of this projectivized bundle.

\begin{thm}[\cite{Demailly95}]
Suppose that $\rank V\ge 2$. Then the quotient $J_k^{\text{reg}}V/\mathbb G_k$ has the structure of a locally trivial bundle over $X$, and there is a holomorphic embedding $J_k^{\text{reg}}V/\mathbb G_k\hookrightarrow X_k$ over $X$, which identifies $J_k^{\text{reg}}V/\mathbb G_k$ with $X_k^{\text{reg}}$, that is the set of point in $X_k$ on the form $f_{[k]}(0)$ for some non singular $k$-jet $f$. In other word $X_k$ is a relative compactification of $J_k^{\text{reg}}V/\mathbb G_k$ over $X$.

Moreover, we have the direct image formula
$$
(\pi_{0,k})_*\mathcal O_{X_k}(m)=\mathcal O(E_{k,m}V^*).
$$
\end{thm}

Next, we are in position to point out the link between the theory of hyperbolicity and invariant jet differentials.

\begin{thm}[\cite{G-G79},\cite{Demailly95}]
Assume that there exist integers $k,m>0$ and an ample line bundle $A\to X$ such that 
$$
H^0(X_k,\mathcal O_{X_k}(m)\otimes\pi_{0,k}^*A^{-1})\simeq H^0(X,E_{k,m}V^*\otimes A^{-1})
$$
has non zero sections $\sigma_1,\dots,\sigma_N$ and let $Z\subset X_k$ be the base locus of these sections. Then every entire holomorphic curve $f\colon\mathbb C\to X$ tangent to $V$ is such that $f_{[k]}(\mathbb C)\subset Z$. In other words, for every global $\mathbb G_k$-invariant differential equation $P$ vanishing on an ample divisor, every entire holomorphic curve $f$ must satisfy the algebraic differential equation $P(f)=0$ (and a similar result is true also for the bundle $E_{k,m}^{GG}T^*_X$).
\end{thm}

\subsection{Algebraic Holomorphic Morse Inequalities}

Let $L\to X$ be a holomorphic line bundle over a compact complex manifold of dimension $n$ and $E\to X$ a holomorphic vector bundle of rank $r$. Suppose that $L$ can be written as the difference of two nef line bundles, say $L=F\otimes G^{-1}$, with $F,G\to X$ numerically effective. Then we have the following asymptotic estimate for the partial alternating sum of the dimension of cohomology groups of $L$ with values in $E$.

\begin{thm}[\cite{Demailly00}]
With the previous notation, we have (strong algebraic holomorphic Morse inequalities):
$$
\sum_{ j=0}^q(-1)^{q-j}h^j(X,L^{\otimes m}\otimes E)\le r\frac{m^n}{n!}\sum_{j=0}^q(-1)^{q-j}\binom nj F^{n-j}\cdot G^j+o(m^n).
$$
In particular \cite{Trapani95}, $L^{\otimes m}\otimes E$ has a global section for $m$ large if $F^n-nF^{n-1}\cdot G>0$.
\end{thm}

\subsection{Schur Powers of a Complex Vector Space}

Here, we just recall the notation and a possible construction of Schur powers of a complex vector space. Let $V$ be a complex vector space of dimension $r$. With every set of nonincreasing $r$-tuples  $(a_1,\dots,a_r)\in\mathbb Z^r$, $a_1\ge a_2\ge\cdots\ge a_r$, one associates, in a functorial way, a collection of vector spaces $\Gamma^{(a_1,\dots,a_r)}V$, which provide the list of all irreducible representations of the linear group Gl$(V)$, up to isomorphism (in fact $(a_1,\dots,a_r)$ is the highest weight of the action of a maximal torus $(\mathbb C^*)^r\subset\text{Gl}(V)$). The Schur functors can be defined in an elementary way as follows. Let 
$$
\mathbb U_r=\left\{\begin{pmatrix} 1 & 0 \\ * & 1 \end{pmatrix} \right\}
$$
be the group of lower triangular unipotent $r\times r$ matrices. If all $a_j$ are nonnegative, one defines
$$
\Gamma^{(a_1,\dots,a_r)}V\subset S^{a_1}\otimes\cdots S^{a_r}V
$$
as the set of polynomials $P(\xi_1,\dots,\xi_r)$ on $(V^*)^r$ which are homogeneous of degree $a_j$ with respect to $\xi_j$ and which are invariant under the left action of $\mathbb U_r$ on $(V^*)^r=\text{Hom}(V,\mathbb C^r)$, i.e.,
$$
P(\xi_1,\dots,\xi_{j-1},\xi_{j}+\xi_k,\xi_{j+1},\dots,\xi_r)=P(\xi_1,\dots,\xi_r),\quad\forall k<j.
$$
We agree that $\Gamma^{(a_1,\dots,a_r)}V=0$ unless $(a_1,\dots,a_r)$ is nonincreasing. As a special case, we recover symmetric and exterior powers
$$
\begin{aligned}
& S^mV=\Gamma^{(m,0,\dots,0)}V, \\
& \bigwedge{}\!\!^k V=\Gamma^{(1,\dots,1,0,\dots,0)}V,\quad\text{with $k$ indices $1$}.
\end{aligned}
$$
The Schur functors satisfy the well-known formula
$$
\Gamma^{(a_1+\ell,\dots,a_r+\ell)}V=\Gamma^{(a_1,\dots,a_r)}V\otimes(\det V)^\ell,
$$
which can be used to define $\Gamma^{(a_1,\dots,a_r)}V$ is any of the $a_j$'s happens to be negative.

\section{Proof of Theorem 1}

First of all, we recall a theorem contained in \cite{B-R90}:

\begin{thm}\label{BR}
Let $X$ be a $n$-dimensional smooth complete intersection in $\mathbb P^N$. Let $T$ be any Young tableau and $t_i$ be the number of cells inside the $i$-th column of T. Set
$$
t:=\sum_{i=1}^{N-n}t_i,\quad\text{$t_i=0$ if $i>\length T$}.
$$
Then, if $t<n$ one has the vanishing
$$
H^0(Y,\Gamma^TT^*_X)=0,
$$
i.e. the smooth complete intersection $X$ has no global $T$-symmetrical tensor forms different from zero if the Young tableau $T$ has less than $\dim X$ cells inside its $\codim X$ front columns.
\end{thm}

In our notation, the irreducible Gl$(T^*_X)$-representation, given for $\lambda_1\ge\lambda_2\ge\dots\ge\lambda_n$ by $\Gamma^{(\lambda_1,\dots,\lambda_n)}T^*_X$, corresponds to $\Gamma^{T_\lambda}T^*_X$ where the tableau $T_\lambda$ is obtained from the partition $\lambda_1+\dots+\lambda_n$. Thus, for example, the tableau with only one row of length $m$ corresponds to $S^mT^*_X$ and the tableau with only one column of depth $k$ corresponds to $\bigwedge^kT^*_X$.

We now need an algebraic lemma.

\begin{lem}
Let $V$ be a complex vector space of dimension $n$ and $\lambda=(\lambda_1,\dots,\lambda_n)$ such that $\lambda_1\ge\lambda_2\ge\dots\ge\lambda_n\ge 0$. Then
$$
\Gamma^\lambda V\otimes S^mV\simeq\bigoplus_\mu\Gamma^\mu V
$$
as Gl$(V)$-representations, the sum being over all $\mu$ whose Young diagram is obtained by adding $m$ boxes to the Young diagram of $\lambda$, with no two in the same column.
\end{lem}

\begin{proof}
This follows immediately by Pieri's formula, see e.g. \cite{F-H91}.
\end{proof}

Note that this implies that between the irreducible Gl$(V)$-representa\-tions of $S^lV\otimes S^m V$, we cannot find terms of type $\Gamma^{(\lambda_1,\dots,\lambda_n)}V$ with $\lambda_i>0$ for $i>2$ (they are all of type $\Gamma^{(l+m-j,j,0,\dots,0)}V$ for $j=0,...,\min\{m,l\}$).

Thus, by induction on the number of factor in the tensor product of symmetric powers,  we easily find:

\begin{prop}
If $k\le n$ then we have a direct sum decomposition into irreducible $\text{Gl}(V)$-represen\-tations
$$
S^{l_1}V\otimes S^{l_2}V\otimes\dots\otimes S^{l_k}V=\bigoplus_\lambda\nu_\lambda\,\Gamma^\lambda V,
$$
where $\nu_\lambda\ne0$ only if $\lambda=(\lambda_1,\dots,\lambda_n)$ is such that $\lambda_i=0$ for $i>k$.
\end{prop}

Now, we can state and prove a slightly more general version (which actually implies immediately Theorem 1) of our first result.

\begin{thm} 
Let $X\subset\mathbb P^{N}$ be a smooth complete intersection. Then 
$$
H^0(X,E_{k,m}^{GG}T^*_X)=0
$$
for all $m\ge 1$ and $1\le k< \dim X/\codim X$.
\end{thm}

\begin{proof}
The bundle $E_{k,m}^{GG}T^*_X$ admits a filtration whose associated graded bundle is given by
$$
\operatorname{Gr}^\bullet E_{k,m}^{GG}T^*_X=\bigoplus_{l_1+2l_2+\dots+kl_k=m}S^{l_1}T^*_X\otimes S^{l_2}T^*_X\otimes\dots\otimes S^{l_k}T^*_X.
$$
The addenda in the direct sum decomposition into irreducible Gl$(T^*_X)$-representation of $\operatorname{Gr}^\bullet E_{k,m}^{GG}T^*_X$ are all of type $\Gamma^{(\lambda_1,\dots,\lambda_n)}T^*_X$ with $\lambda_i=0$ for $i>k$. This means that in the statement of Theorem \ref{BR}, each $t_i$ is less than or equal to $k$. Then, if $\dim X=n$, we have
$$
\begin{aligned}
t & =\sum_{i=1}^{N-n}t_i \le \sum_{i=1}^{N-n} k \\
& = k(N-n) < \frac n{N-n}(N-n)=n.
\end{aligned}
$$ 
Thus, by Theorem \ref{BR}, we have vanishing of global section of each graded piece. 
We now only need to link the vanishing of the cohomology of a filtered vector bundle to the vanishing of his graded bundle. This is done in the next lemma and the proof is achieved.
\end{proof}

\begin{lem}
Let $E\to X$ be a holomorphic filtered vector bundle with filtered pieces $\{0\}=E_r\subset\cdots\subset E_{p+1}\subset E_p\subset\cdots\subset E_0=E$. If $H^q(X,\operatorname{Gr}^\bullet E)=0$ then $H^q(X,E)=0$.
\end{lem}

\begin{proof}
Consider the short exact sequence
$$
0\to \operatorname{Gr}^p E\to E/E_{p+1}\to E/E_p\to 0
$$
and the associated long exact sequence in cohomology
$$
\cdots\to \underbrace{H^q(X,\operatorname{Gr}^pE)}_{=0}\to H^q(X,E/E_{p+1})\to H^q(X, E/E_p)\to\cdots.
$$
For $p=1$, we get $H^q(X, E/E_1)=H^q(X, \operatorname{Gr}^0E)=0$ by hypothesis. Thus, $H^q(X, E/E_2)=0$ and, by induction on $p$, we find $H^q(X, E/E_p)=0$. Therefore, for $p=r$ we get the desired result.
\end{proof}

\section{Proof of Theorem 2}

In this section, we prove our Theorem 2 by means of the algebraic version of holomorphic Morse inequalities, performed on a particular subbundle of the invariant jet differentials. For the sake of simplicity, we shall reproduce here only the proof for three jets on threefolds in $\mathbb P^4$, the other cases being completely analogous.

\subsection{Chern Classes Computations}

Denote by $c_\bullet(E)$ the total Chern class of a vector bundle $E$. The short exact sequences (\ref{ses1}) and (\ref{ses2}) give us, for each $k>0$, the following formulae:
$$
c_\bullet(V_k)=c_\bullet(T_{X_k/X_{k-1}}) c_\bullet(\mathcal O_{X_k}(-1))
$$
and
$$
c_\bullet(\pi_k^*V_{k-1}\otimes\mathcal O_{X_k}(1))=c_\bullet(T_{X_k/X_{k-1}}),
$$
so that
\begin{equation}\label{chern1}
c_\bullet(V_k)=c_\bullet(\mathcal O_{X_k}(-1))c_\bullet(\pi_k^*V_{k-1}\otimes\mathcal O_{X_k}(1)).
\end{equation}
Let us call $u_j=c_1(\mathcal O_{X_j}(1))$ and $c_l^{[j]}=c_l(V_j)$. With this notation, (\ref{chern1}) becomes
\begin{equation}\label{chern2}
1+c_1^{[k]}+\cdots+c_r^{[k]}=(1-u_k)\sum_{0\le j\le r}\pi_k^*c_j^{[k-1]}(1+u_k)^{r-j}.
\end{equation}
Since $X_j$ is the projectivized bundle of line of $V_{j-1}$, we also have the polynomial relations
\begin{equation}\label{chern3}
u_j^r+\pi_j^*c_1^{[j-1]}\cdot u_j^{r-1}+\cdots+\pi_j^*c_{r-1}^{[j-1]}\cdot u_j+\pi_j^*c_{r}^{[j-1]}=0,\quad 1\le j\le k.
\end{equation}

\begin{prop}\label{relchern}
Let $\rank V=3$. Then, we have the following relations for Chern classes:
\begin{equation}\label{chernV3}
\begin{aligned}
& c_1^{[k]}=\pi_k^*c_1^{[k-1]}+2u_k, \\
& c_2^{[k]}=\pi_k^*c_2^{[k-1]}+\pi_k^*c_1^{[k-1]}\cdot u_k, \\
& c_3^{[k]}=\pi_k^*c_3^{[k-1]}-\pi_k^*c_1^{[k-1]}\cdot u_k^2-2u_k^3, \\
& u_k^3+\pi_k^*c_1^{[k-1]}\cdot u_j^{2}+\pi_k^*c_{2}^{[k-1]}\cdot u_k+\pi_k^*c_{3}^{[k-1]}=0.
\end{aligned}
\end{equation}
\end{prop}

\begin{proof}
This is just a straightforward computation using identity (\ref{chern2}).
\end{proof}

Next, let $X\subset\mathbb P^{n+1}$ be a smooth hypersurface of degree $\deg X=d$. Then, we have a short exact sequence
$$
0\to T_X\to T_{\mathbb P^{n+1}}|_X\to \mathcal O_X(d)\to 0;
$$
so we have the following relation for the total Chern class of $X$
$$
(1+h)^{n+2}=(1+d\,h)c_\bullet(X),
$$
where $h=c_1(\mathcal O_{\mathbb P^{n+1}}(1))$ and $(1+h)^{n+2}$ is the total Chern class of $\mathbb P^{n+1}$. Thus, an easy computation yields:

\begin{prop}\label{degchern}
Let $X\subset\mathbb P^{4}$ be a smooth hypersurface of degree $\deg X=d$. Then, the Chern classes of $X$ are given (in term of the hyperplane divisor) by
\begin{equation}\label{chernX3}
\begin{aligned}
& c_1(X)=-h(d-5), \\
& c_2(X)=h^2(d^2-5d+10), \\
& c_3(X)=-h^3(d^3-5d^2+10d-10) \\
\end{aligned}
\end{equation}
and $h^3=d$.
\end{prop}

\subsection{Choice of the Appropriate Subbundle}

By definition, there is a canonical injection $\mathcal O_{X_k}(-1)\hookrightarrow\pi_k^*V_{k-1}$ and a composition with the differential of the projection $(\pi_k)_*$ yields, for all $k\ge 2$, a canonical line bundle morphism
$$
\mathcal O_{X_k}(-1)\hookrightarrow\pi_k^*V_{k-1}\to\pi_k^*\mathcal O_{X_{k-1}}(-1),
$$
which admits precisely $D_k\overset{\text{def}}=P(T_{X_{k-1}/X_{k-2}})\subset P(V_{k-1})=X_k$ as its zero divisor. Hence, we find
\begin{equation}\label{mor}
\mathcal O_{X_k}(1)=\pi_k^*\mathcal O_{X_{k-1}}(1)\otimes\mathcal O(D_k).
\end{equation}
Now, for $\bold a=(a_1,\dots,a_k)\in\mathbb Z^k$, define a line bundle $\mathcal O_{X_k}(\bold a)$ on $X_k$ as
$$
\mathcal O_{X_k}(\bold a)=\pi_{1,k}^*\mathcal O_{X_1}(a_1)\otimes\pi_{2,k}^*\mathcal O_{X_2}(a_2)\otimes\cdots\otimes\mathcal O_{X_k}(a_k).
$$
By (\ref{mor}), we have
$$
\pi_{j,k}^*\mathcal O_{X_j}(1)=\mathcal O_{X_k}(1)\otimes\mathcal O_{X_k}(-\pi_{j+1,k}^*D_{j+1}-\cdots-D_k),
$$
thus by putting $D_j^\star=\pi_{j+1,k}^*D_{j+1}$ for $j=1,\dots,k-1$ and $D_k^\star=0$, we have an identity
$$
\begin{aligned}
& \mathcal O_{X_k}(\bold a)=\mathcal O_{X_k}(b_k)\otimes\mathcal O_{X_k}(-\bold b\cdot D^\star),\quad\text{where} \\
& \bold b=(b_1,\dots,b_k)\in\mathbb Z^k,\quad b_j=a_1+\cdots+a_j, \\
& \bold b\cdot D^\star=\sum_{j=1}^{k-1}b_j\,\pi_{j+1,k}^*D_{j+1}.
\end{aligned}
$$
In particular, if $\bold b\in\mathbb N^k$, that is if $a_1+\cdots+a_j\ge 0$, we get a morphism
$$
\mathcal O_{X_k}(\bold a)=\mathcal O_{X_k}(b_k)\otimes\mathcal O_{X_k}(-\bold b\cdot D^\star)\to\mathcal O_{X_k}(b_k).
$$
We then have the following:

\begin{prop}[\cite{Demailly95}]
Let $\bold a=(a_1,\dots,a_k)\in\mathbb N^k$ and $m=a_1+\cdots+a_k$. 
\begin{itemize}
\item We have the direct image formula
$$
(\pi_{0,k})_*\mathcal O_{X_k}(\bold a)\simeq\mathcal O(\overline F^{\bold a}E_{k,m}V^*)\subset\mathcal O(E_{k,m}V^*)
$$
where $\overline F^{\bold a}E_{k,m}V^*$ is the subbundle of polynomials $Q(f',f'',\dots,f^{(k)})\in E_{k,m}V^*$ involving only monomials $(f^{(\bullet)})^\ell$ such that
$$
\ell_{s+1}+2\ell_{s+2}+\cdots+(k-s)\ell_{k}\le a_{s+1}+\cdots+a_k
$$
for all $s=0,\dots,k-1$.
\item If $a_1\ge 3a_2,\dots,a_{k-2}\ge 3a_{k-1}$ and $a_{k-1}\ge 2a_k> 0$, the line bundle $\mathcal O_{X_k}(\bold a)$ is relatively nef over $X$.
\end{itemize}
In particular, the line bundle $\mathcal L_k(X)\overset{\text{def}}=\mathcal O_{X_k}(2\cdot 3^{k-2},2\cdot 3^{k-3},\dots,6,2,1)$ is relatively nef over $X$ and its direct image on $X$ is a subbundle of the bundle of invariant jet differentials of order $k$ and weighted degree $3^{k-1}$.
\end{prop}

In the case of projective hypersurface, we obtain the following expression of $\mathcal L_k$ as the difference of two nef line bundles.

\begin{lem}
Let $X\subset\mathbb P^{n+1}$ be a projective hypersurface. Then $\mathcal L_k(X)\otimes\pi_{0,k}^*\mathcal O_X(l)$ is nef if $l\ge 2\cdot 3^{k-1}$. In particular, 
$$
\mathcal L_k(X)=\mathcal F_k(X)\otimes\mathcal G_k(X)^{-1},
$$
where $\mathcal F_k(X)=\mathcal L_k(X)\otimes\pi_{0,k}^*\mathcal O_X(2\cdot 3^{k-1})$ and $\mathcal G_k(X)=\pi_{0,k}^*\mathcal O_X(2\cdot 3^{k-1})$ are nef. 
\end{lem}

\begin{proof}
Of course, as a pull-back of an ample line bundle, 
$$
\mathcal G_k(X)=\pi_{0,k}^*\mathcal O_X(2\cdot 3^{k-1})
$$ 
is nef. It is well known that the cotangent space of the projective space twisted by $\mathcal O(2)$ is globally generated. Hence, $T^*_X\otimes\mathcal O_X(2)$ is globally generated as a quotient of $T^*_{\mathbb P^{n+1}}|_X\otimes\mathcal O_X(2)$, so that $\mathcal O_{X_1}(1)\otimes\pi_{0,1}^*\mathcal O_X(2)=\mathcal O_{\mathbb P(T^*_X\otimes\mathcal O_X(2))}(1)$ is nef. 

Next, we construct by induction on $k$, a nef line bundle $A_{k}\to X_{k}$ such that $\mathcal O_{X_{k+1}}(1)\otimes\pi_k^*A_{k}$ is nef. By definition, this is equivalent to say that the vector bundle $V_{k}^*\otimes A_{k}$ is nef. By what we have just seen, we can take $A_0=\mathcal O_X(2)$ on $X_0=X$. Suppose $A_0,\dots,A_{k-1}$ as been constructed. As an extension of nef vector bundles is nef, dualizing the short exact sequence (\ref{ses1}) we find
$$
0\to\mathcal O_{X_k}(1)\to V_k^*\to T_{X_k/X_{k-1}}^*\to 0,
$$ 
and so we see, twisting by $A_k$, that it suffices to select $A_k$ in such a way that both $\mathcal O_{X_k}(1)\otimes A_k$ and $T_{X_k/X_{k-1}}^*\otimes A_k$ are nef. To this aim, considering the second wedge power of the central term in (\ref{ses2}), we get an injection
$$
0\to T_{X_k/X_{k-1}}\to\bigwedge{}\!\!^2(\pi_k^* V_{k-1}\otimes\mathcal O_{X_k}(1))
$$
and so dualizing and twisting by $\mathcal O_{X_{k}}(2)\otimes\pi_k^* A_{k-1}^{\otimes 2}$, we find a surjection
$$
\pi_k^*\bigwedge{}\!\!^2(V_{k-1}^*\otimes A_{k-1})\to T^*_{X_k/X_{k-1}}\otimes\mathcal O_{X_k}(2)\otimes\pi_k^* A_{k-1}^{\otimes 2}\to 0.
$$
By induction hypothesis, $V_{k-1}^*\otimes A_{k-1}$ is nef so the quotient $T^*_{X_k/X_{k-1}}\otimes\mathcal O_{X_k}(2)\otimes\pi_k^* A_{k-1}^{\otimes 2}$ is nef, too.
In order to have the nefness of both $\mathcal O_{X_k}(1)\otimes A_k$ and $T_{X_k/X_{k-1}}^*\otimes A_k$, it is enough to choose $A_k$ in such a way that $A_k\otimes\pi_k^* A_{k-1}^*$ and $A_k\otimes\mathcal O_{X_k}(-2)\otimes\pi_k^*{A_{k-1}^*}^{\otimes 2}$ are both nef: therefore we set 
$$
A_k=\mathcal O_{X_k}(2)\otimes\pi_k^* A_{k-1}^{\otimes 3}=\bigl (\mathcal O_{X_k}(1)\otimes\pi_k^* A_{k-1}\bigr)^{\otimes 2}\otimes\pi_k^* A_{k-1},
$$
which, as a product of nef line bundles, is nef and satisfies the two conditions above.
This gives $A_k$ inductively, and the resulting formula for $\mathcal O_{X_k}(1)\otimes\pi_k^*A_{k-1}$ is
$$
\begin{aligned}
\mathcal O_{X_k}(1)\otimes\pi_k^*A_{k-1} &= \mathcal L_{k}(X)\otimes\pi_{0,k}^*\mathcal O_X(2\cdot(1+2+\cdots+2\cdot 3^{k-2})) \\
& =   \mathcal L_{k}(X)\otimes\pi_{0,k}^*\mathcal O_X(2\cdot 3^{k-1}).
\end{aligned}
$$
The lemma is proved.
\end{proof}

\subsection{End of the Proof of Theorem 2} 
We just apply the algebraic version of holomorphic Morse inequalities to the line bundle $\mathcal L_3(X)$ for $X$ a smooth hypersurface in $\mathbb P^4$ (and, for the other cases, to the line bundle $\mathcal L_n(X)$ for $X$ a smooth hypersurface in $\mathbb P^{n+1}$, $n=2,\dots,5$).

Then, we have to compute $\mathcal F_3(X)^9-9\,\mathcal F_3(X)^8\cdot\mathcal G_3(X)$ which is given in terms of Chern classes by
$$
\begin{aligned}
& (u_3+2\,\pi_{2,3}^*u_2+6\,\pi_{1,3}^*u_1+18\,\pi_{0,3}^*h)^9 \\
&\quad -9(u_3+2\,\pi_{2,3}^*u_2+6\,\pi_{1,3}^*u_1+18\,\pi_{0,3}^*h)^8\cdot 18\,\pi_{0,3}^*h.
\end{aligned}
$$
By using recursively Proposition \ref{relchern} we can express this quantity in terms of Chern classes of $X$ and the hyperplane class (the computation is made with GP/PARI CALCULATOR Version 2.3.2):
$$
\begin{aligned}
-3421377792\,h^3 + 676045440\,c_1(X)\cdot h^2 -7494966\,c_1(X)^3  \\
+ 10997352\,c_2(X)\cdot c_1(X) - 3835548\,c_3(X),
\end{aligned}
$$
and, by Proposition \ref{degchern}, we obtain
$$
\begin{aligned}
\mathcal F_3(X)^9-9\,&\mathcal F_3(X)^8\cdot\mathcal G_3(X) \\ 
& = 333162\,d^4 - 21628710\,d^3 - 460474830\,d^2 - 466509222\,d,
\end{aligned}
$$
which is positive if $d=\deg(X)\ge 82$.

\begin{rem}
Even if we know \cite{Rousseau06b} that the line bundle $\mathcal O_{X_3}(1)$ is big in the case of smooth hypersurfaces in projective $4$-space for $\deg(X)\ge 97$, to get the result with these techniques we are obliged to deal with the line bundle $\mathcal L_3(X)$. In fact, the algebraic version of holomorphic Morse inequalities gives, if we merely utilise $\mathcal O_{X_3}(1)$, a negative lower bound (and similarly in higher dimension). The reason is that $\mathcal O_{X_k}(1)$ is always relatively big but never relatively nef (if $k\ge 2$) so that with $\mathcal O_{X_k}(1)$, holomorphic Morse inequalities take into account too many \lq\lq relatively negative terms\rq\rq{}.
\end{rem}

\begin{rem}
With a slight theoretical refinement of these techniques, we are able to show the existence of global invariant jet differentials vanishing on an ample divisor even for the dimension of the smooth hypersurface less than or equal to $8$. Unfortunately, the trick we use to increase the dimension leads to the loss of the effectivity:  we can show the existence of a lower bound for the degree depending only on the dimension of the hypersurface, but we are unable to compute it.

However, we strongly feel that in a forthcoming future, with a further improvement in the comprehension of the cohomology algebra of $X_k$, we will be able to show (without explicit effectivity) the existence of global invariant jet differentials vanishing on an ample divisor for smooth hypersurfaces of any dimension.
\end{rem}


\section{Computer Calculations}
Here we reproduce the code we utilised to perform all computations with GP/PARI CALCULATOR Version 2.3.2.
\newline
\newline
\texttt{
/* scratch variable */ 
\newline 
X
\newline \newline 
/* main formal variables */
\newline 
c=[c1,c2,c3,c4,c5,c6,c7,c8,c9]  /* Chern classes of V on X */
\newline 
u=[u1,u2,u3,u4,u5,u6,u7,u8,u9]  /* Chern classes of OXk(1) */
\newline 
v=[v1,v2,v3,v4,v5,v6,v7,v8,v9]  /* Chern classes of Vk on Xk */
\newline 
w=[w1,w2,w3,w4,w5,w6,w7,w8,w9]  /* formal variables */
\newline 
e=[0,0,0,0,0,0,0,0,0]  /* empty array for Chern classes of hypersurfaces */
\newline 
q=[0,0,0,0,0,0,0,0,0]  /* empty array for Chern equations */
\newline 
 \newline /* main */
\newline Calcul(dim,order)=
\newline \{
\newline local(j,n,N);
\newline n=dim;
\newline r=dim;
\newline k=order;
\newline N=n+k*(r-1);
\newline H(n+1);
\newline Chern();
\newline B=2*h*3\^{}(k-1);
\newline A=B+u[k];
\newline for(j=1,k-1,A=A+2*3\^{}(k-j-1)*u[j]);
\newline R=Reduc((A-N*B)*A\^{}(N-1));
\newline C=Integ(R);
\newline print("Calculation for order ", k, " jets on a ", n, "-fold");
\newline print("Line bundle A= ", A);
\newline print("Line bundle B= ", B);
\newline print("Chern class of A\^{}", N, "-", N, "*A\^{}", N-1, "*B :");
\newline print(C);
\newline E=Eval(C);
\newline print("Evaluation for degree d hypersurface in P\^{}", n+1, " :");
\newline print(E)
\newline 	\}
\newline
\newline
/* compute Chern relations */
\newline Chern()=  
\newline \{
\newline local(j,s,t);
\newline q[1]=X\^{}r; for(j=1,r,q[1]=q[1]+c[j]*X\^{}(r-j));
\newline for(s=1,r,v[s]=c[s]);
\newline for(s=r+1,9,v[s]=0);
\newline for(t=1,k-1,\
  \newline for(s=1,r,w[s]=v[s]+(binomial(r,s)-binomial(r,s-1))*u[t]\^{}s;
  \newline   for(j=1,s-1,w[s]=w[s]+\
  \newline     (binomial(r-j,s-j)-binomial(r-j,s-j-1))*v[j]*u[t]\^{}(s-j)));
\newline  for(s=1,r,v[s]=w[s]);
  \newline q[t+1]=X\^{}r; for(j=1,r,q[t+1]=q[t+1]+v[j]*X\^{}(r-j)))
\newline \}
\newline 
\newline /* reduction to Chern classes of X */
\newline
Reduc(p)= 
\newline \{
\newline local(j,a);
\newline a=p;
\newline for(j=0,k-1,\
\newline   a=subst(a,u[k-j],X);
\newline   a=subst(lift(Mod(a,q[k-j])),X,u[k-j]));
\newline a
\newline  \}
\newline 
\newline /* integration along fibers */
\newline Integ(p)= 
\newline \{
\newline local(j,a);
\newline a=p;
\newline for(j=0,k-1,\
 \newline  a=polcoeff(a,r-1,u[k-j]));
\newline a
\newline \}
\newline
\newline
/* compute Chern classes of hypersurface of degree d in P\^{}n */
\newline H(n)= 
\newline  \{
\newline local(j,s);
\newline for(s=1,n-1,\
  \newline e[s]=binomial(n+1,s);
 \newline  for(j=1,s,e[s]=e[s]+(-d)\^{}j*binomial(n+1,s-j)))
\newline  \}
\newline
\newline
 /* evaluation in terms of the degree */
\newline Eval(p)=
\newline \{
\newline local(a,s);
\newline a=p;
\newline for(s=1,r,a=subst(a,c[s],e[s]));
\newline subst(a,h,1)*d
\newline \}
}

\end{document}